\documentclass{article}

\newtheorem{theorem}{Theorem}[section]
\newtheorem{lemma}{Lemma}[section]
\newtheorem{remark}{Remark}

\title{Boundary feedback stabilization of
the isothermal Euler-equations with uncertain boundary data
}
\author{Martin Gugat
\thanks{
Lehrstuhl Angewandte Mathematik 2,
Department Mathematik,
Friedrich-Alexander Universit\"at Erlangen-N\"urnberg (FAU),
Cauerstr. 11, 91058 Erlangen, Germany,
(email: {martin.gugat@fau.de})
}
\and  R\"udiger Schultz
\thanks{
Universit\"at Duisburg-Essen, Fakult\"at f\"ur Mathematik,
Thea-Leymann-Str. 9, 45127 Essen, Germany}
}
\date{\today}

\begin{document}

\maketitle

\begin{abstract}
In a gas transport system, the
customer behavior is uncertain.
Motivated by this situation, we consider
a boundary stabilization problem
for the flow through  a gas pipeline,
where the outflow at one
end of the pipe
is uncertain.
The  control action is located
at the
other end of the pipe.
The feedback law is
a classical
Neumann velocity feedback with a feedback parameter $k>0$.

We show that
as long as the
$H^1$-norm of the function that describes the noise in the customer's
behavior decays exponentially with
a rate that is sufficiently large,
the velocity  of the gas can be stabilized exponentially
fast in the sense that a suitably chosen Lyapunov function decays exponentially.
For the exponential stability it is sufficient
that the feedback parameter $k$ is sufficiently large
and the stationary state to which the system is stabilized is
sufficiently small.
The stability result is local, that is it holds for initial
states that are sufficiently close to the stationary state.

This result is an example for the exponential boundary feedback stabilization of
a quasilinear hyperbolic system with uncertain boundary data.
The analysis is based upon  the choice of a suitably Lyapunov function.
The decay of this Lyapunov function implies that
also the $L^2$-norm of the
difference of the system state and the stationary state decays exponentially.

%
\end{abstract}



\section{Introduction}
There are many studies on boundary feedback stabilization
where exponential stability is shown
for a family initial states under
the assumption
that the action at the boundary is
certain.
In this paper we consider a different situation,
where not only the initial state,
but also the behavior at a part of
the boundary of the system is uncertain.
We  show that
if the uncertain boundary action
converges sufficiently fast to a
stationary state,
during this process  a suitably chosen Lyapunov function
decays  exponentially.
As stated in \cite{prieur}, where semilinear parabolic systems are considered,
{\em Lyapunov function based techniques are central in the study of partial differential equations} (pdes).
In particular, in  \cite{prieur} input--to-state stable (ISS) Lyapunov functions for pdes
with disturbances are considered.
In this paper we consider a system with a quasilinear hyperbolic pde with
uncertainties in the boundary data.
The concept of input--to--state stability is discussed in detail in
\cite{eduardosontag}.

The transient behavior of pipeline gas flow
can be modeled very precisely by a
quasilinear
hyperbolic partial differential equation.
%
%
%
%
Often the consumer behavior in gas transport networks is uncertain.
Therefore
we consider a problem of  boundary stabilization with uncertain
boundary data
\[b(t,\, \omega)\]
that depend on an uncertain decision $\omega \in \Omega$
that models the  uncertainty in the customer's behavior.
This customer behavior can be considered
as a noise in the boundary data.
We assume that the consumer behavior
approaches some desired state in $H^1$,
that is that the noise has the following structure:

For all $\omega\in \Omega$ the function $b(\cdot,\,\omega)$ is twice continuously differentiable and
there exists  numbers
$T_{period}>0$  and $T>T_{period}
$, such that
 for all $t\in (
T_{period}, \,T)$
\begin{equation}
\label{bassumption}
\int_{t - T_{period}}^t
\left|b(\tau,\, \omega)\right |^2
+
\left|b_t(\tau,\, \omega)\right |^2 d\tau
\leq
\,
C_\nu
\,
\exp\left(- \nu\, t \right)
\end{equation}
for some
$\nu > 0$ and $C_\nu>0$.
%
%
This means that after some variations the
noise in the customer behavior
decays exponentially fast to zero.

We consider the
questions:
How does the uncertainty in the demand influence the stabilization of
the system?
Is there a boundary feedback law that leads to
exponential decay of the difference of
the system state and a desired stationary state
in spite of the uncertain boundary data?

For this purpose we consider a boundary feedback law that
for the case without noise, that is
for $b(t,\, \omega)=0$ for all $t\geq 0$, stabilizes the system exponentially fast.


Output-feedback stabilization of stochastic nonlinear systems
driven by noise of unknown covariance has been studied in  \cite{krstic}
where a controller is presented that guarantees regulation
to the desired state with probability one.
Our result is of a different type.
We present a Lyapunov function
that decays exponentially under
suitable smallness assumptions for
the variations in the consumer behavior.
Our Lyapunov function for the case with uncertainty
is a time average of
a strict Lyapunov function
for the case without uncertainty.

In \cite{DGL1} a
 strict $H^1$-Lyapunov function and feedback stabilization for the isothermal Euler equations with friction
 have been studied without uncertainty.
 The stabilization of
 the gas flow in pipeline networks has been considered in \cite{dick}.
 The novelty in the present contribution is that we include uncertainty in
 the boundary data in our analysis.
The isothermal Euler equations are equivalent to a quasilinear wave equation for
the gas velocity.
Therefore for the stabilization of our system we consider
the same feedback law that is used for the stabilization for the linear wave equation
that has been studied for example in \cite{ko:rapid}.
For the linear wave equation, this feedback law
stabilizes the system as long as the feedback parameter
has the right sign.

In our system, due to its nonlinearity, the stationary states are not constant
and blow up after a finite  critical length.
Due to the nonlinearity of the system, in our analysis we have to assume that the stationary
states are sufficiently small and we can only stabilize
the system locally around the stationary state.
To guarantee the exponential decay of the Lyapunov function,
we have to assume that the feedback parameter is sufficiently large.
%
Apart from the obvious restriction that
the length of the pipe is less than the critical length,
in this paper we
do not impose any additional restrictions on the length of
the pipe.
This is in contrast to the earlier contributions
\cite{guherty} and \cite{GLTW} that are only applicable if
the lengths of the pipes are sufficiently small.


\section{The model for the pipeline flow}

Let a finite time $T>0$ be given.
The system dynamics for the  flow of ideal gas in a single pipe can be modeled
by
the isothermal Euler
equations (see \cite{BHK1},\cite{BHK2},\cite{DGL1}):
\begin{eqnarray}
\label{2.1}
& &\rho_t+q_x=0,\\
\label{2.2}
& &q_t+\Big(\frac{q^2}\rho+a^2\rho\Big)_x=-\frac{1}{2}\,\theta \,
\frac{q|q|}{\rho}
\end{eqnarray}
where $\rho=\rho(t,x)>0$ is the density of the gas, $q=q(t,x)$ is
the mass flux,
$\theta=\frac{f_g}{\delta}$ where
 the constant $f_g \geq  0$ is a friction factor and
$\delta>0$ is the diameter of the pipe.
The constant $a > 0$ is the speed of sound
in the gas. We consider the equations on the domain
$\Omega:=[0,T]\times[0,L]$ where $L>0$ and $T>0$ are given.
Equation (\ref{2.1}) states the conservation of mass and equation
(\ref{2.2}) is the balance of  momentum.
Define the velocity   $\tilde u$  of the gas flow as
\begin{equation}
\label{2.35}
\tilde u=\frac{q}\rho.
\end{equation}
In this paper, we consider subsonic  positive gas flow,
that is we assume that
\begin{equation}
\label{2.3}
0<
\tilde u
 < a.
\end{equation}
Note that in the operation of the gas pipelines,
there are strict upper bounds for the velocities
in order to avoid noise pollution
by pipeline vibrations that can be generated by the flow,
see \cite{zou}.
%
%
%
%
For sufficiently regular states, $\tilde u$ satisfies the
quasilinear wave equation (see \cite{gubook})
\begin{equation}
\label{2.21}
\tilde u_{tt} + 2\, \tilde u \, \tilde u_{tx}-(a^2-\tilde u^2)
\, \tilde
u_{xx}=\tilde F(\tilde u,\tilde u_x,\tilde u_t).
\end{equation}
The lower order term
is
\begin{equation}\label{2.22}
\tilde F(\tilde u,\tilde u_x,\tilde u_t)
=
- 2 \, \tilde u_t \, \tilde u_x
-2 \, \tilde u \, \tilde u_x^2
-\frac{3}{2}\, \theta \, \tilde u  \, |\tilde u| \, \tilde u_x
-\theta \, |\tilde u| \, \tilde u_t.
\end{equation}

We consider stationary states $\bar u$ of (\ref{2.21})
that  are
solutions of the ordinary differential equation
\[
\bar u_x=\frac{\theta}{2}\;
\frac{1}{(a^2 - \bar u^2)} \; \, |\bar u|\, \bar u^2.
\]
Obviously in the subsonic case
these solutions are strictly increasing.
These  stationary states  correspond to  stationary states
of the system (\ref{2.1}), (\ref{2.2})
and
are discussed in detail in \cite{nhm}.
In fact with $\sigma \in \{-1,\, 1\}$
that determines the direction of the flow
and a real constant $c_1< -1 $ the stationary
states  have the form
\begin{equation}
\label{stationaryrepresentation}
\bar u(x) = \frac{ \sigma\, a}
{ \sqrt
{
- W_{-1} \left( -\exp\left(\sigma
\,
\theta \, x + c_1  \right)
 \right)
 }
 }
\end{equation}
where $W_{-1}$ is
the Lambert-W   function (see \cite{LA}, \cite{CGH}),
that is the inverse function of
$x\mapsto x \exp(x)$ for $x\leq -1$.
Thus $W_{-1}$ is defined  on $(-\frac{1}{\rm e},  0)$.

For positive gas flow we have $\sigma=1$.
The representation (\ref{stationaryrepresentation})
implies that for $\sigma =1$, the solution $\bar u(x)$ exists only
for $x\leq L_{crit}$ with a critical length $L_{crit}$
and at $L_{crit}$, the flow becomes sonic that is
$\bar u( L_{crit} ) =a$ and the derivative blows up.


To stabilize the system governed by
the quasilinear wave equation (\ref{2.21}) locally around a
desired
stationary state $\bar u(x)$,
at $x=0$ we use the Neumann boundary feedback law
\begin{equation}
\label{0rb}
\tilde u_x(t,\,0)
 =
 \bar u_x(0)+k \, \tilde u_t(t,\,0)
\end{equation}
with a  feedback parameter
$k \in (0,\infty)$.
At $x=L$, the Dirichlet boundary condition is
\begin{equation}
\label{lrb}
\tilde u(t,\,L)
 =
   b(t,\, \omega) + \bar u(L).
\end{equation}
The feedback law
(\ref{0rb})
 is similar to the feedback law for
the stabilization of the linear wave equation
that is studied  for example in \cite{ko:rapid}, \cite{komzu} and \cite{kos}.
%
%
Define
\begin{equation}
u = \tilde u -\bar u
\end{equation}
that is $u$ is the difference of the velocity and the stationary velocity.
Our system stated in terms of $u$ is
\begin{equation}
\label{sys}
\left\{
\begin{array}{l}
u(0,\, x) = \varphi(x),\; x\in [0,\, L]
\\
u_t(0,\, x) = \psi(x),\;  x\in [0,\, L]
\\
u_{tt}+2\,(\bar u+u)\,u_{tx}-\Big(a^2-(\bar u+u)^2\Big)\,u_{xx}=
F(x,u,u_x,u_t)\,\; \mbox{\rm on }\; [0,\,T]\times [0,\, L]
\\
u_x(t,\, 0) =  k\, u_t(t,\, 0),\; t \in [0,\, T]
\\u(t,\, L) = b(t,\omega),\; t \in [0,\, T]
\end{array}
\right.
\end{equation}
where
 $F:=F(x,u,u_x,u_t)$ satisfies
\begin{eqnarray}
\label{Fdefinition}
F & = &
 \tilde F(u + \bar u,\, u_x + \bar u_x\,,u_t)
 -
 \frac{ a^2 - (\bar u + u)^2}{a^2 - \bar u^2} \, \tilde F( \bar u,\,  \bar u_x\,,0).
\end{eqnarray}

\section{Well--posedness}

In \cite{LBW}  a result about
semi-global $C^2$--solutions of quasilinear wave equations is proved
that we can apply to show the well--posedness of (\ref{sys}):
\begin{theorem}
\label{wellposedness}
Let $L>0$, $a>0$,  $k>0$
and a subsonic stationary state
$\bar u\in C^1([0,L])$ be given such that
$\bar u(x) \in (0,\; a)$. 
Choose $T>0$ arbitrarily large.

There exist  constants
$\varepsilon_0(T)>0$
and $C_T>0$,
such that if the
initial data $(u(0,\,x),\, u_t(0,\, x)) = (\varphi(x),\, \psi(x))\in C^2([0,L])\times C^1([0,L])$
and $b(\cdot,\,\omega)$
satisfy
\begin{equation}
\label{epsilon0bedingung}
\max\left\{  \|\varphi(x)\|_{ C^2([0,L])},\;
   \|\psi(x))\|_{ C^1([0,L])},\; \|b(\cdot,\,\omega)\|_{ C^2([0,T])}
   \right\}
   \leq \varepsilon_0(T)
\end{equation}
and the $C^2$-compatibility conditions are satisfied at the points
$(t,x)=(0,0)$ and  $(0,L)$,
then the initial-boundary problem (\ref{sys})
has a unique
solution $u(t,x)\in
C^2([0,T]\times [0,L])$.
%
Moreover the following {\bf a priori estimate} holds:
\begin{equation}
\label{apriori}
\| u\|_{ C^2([0,T]\times [0,L])} \leq C_T \max\left\{  \|\varphi(x)\|_{ C^2([0,L])},\;
   \|\psi(x))\|_{ C^1([0,L])},\; \|b(\cdot,\,\omega)\|_{ C^2([0,T])}
   \right\}.
   \end{equation}
\end{theorem}

\section{Exponential Decay}

In order to show the exponential decay
for the closed--loop system (\ref{sys})  with the quasilinear wave equation, we define
\begin{equation}
\label{e1definition}
E_1(t)
=
\int_0^L
k\, \Big[\left(a^2-(\bar u+u)^2 \right)\, u_x^2+u_t^2\Big]-2
\exp\left(- \frac{x}{L}\right)\,
\Big[(\bar u+u)\,u_x^2+u_t \, u_x\Big]\,dx.
\end{equation}
Note that the first term of $E_1(t)$ is
similar to the classical energy
\[E_{classic}(t) = k\, \int_0^L\, a^2 \, u_x^2 + u_t^2\,dx.\]
The
characteristic curves
for the
quasilinear wave  equation in (\ref{sys})
have
the non-constant slopes $\lambda_+ =\bar u + u + a$
and  $\lambda_- =\bar u + u - a$,
so for the product of the eigenvalues we have
$\left|\lambda_+ \, \lambda_-\right| = a^2  -(\bar u+u)^2$.
For the linear wave equation, this corresponds to
the constant $\left|a\,(-a)\right|$.
Therefore the constant $a^2$ in the definition of $E_{classic}(t)$
is replaced by the function $\left|\lambda_+ \, \lambda_-\right|$
 in the definition of $E_1(t)$.
The second term is added in the integral in $E_1(t)$
since in the quasilinear wave equation in (\ref{sys}),
the mixed partial derivative $u_{tx}$ appears.
The corresponding coefficient
is given by the sum of the eigenvalues
$\lambda_+ +\lambda_- = 2 (\bar u + u)$.
For the linear wave equation, this corresponds to
the constant $a +  (-a)=0$, therefore
the mixed partial derivative does not appear.
Exponential weights
like $\exp\left(- \frac{x}{L}\right)$
 have already been used  
to construct strict Lyapunov functions for hyperbolic pdes, see
for example \cite{bastin}, \cite{Coronbook}.

The following lemma shows that $\sqrt{E_1}$
is equivalent to the $L^2$--norm of $(u_t,\, u_x)$.
\begin{lemma}
\label{fragezeichen}
Assume that
\begin{equation}
\label{110816}
0 \leq \bar u + u \leq \frac{a}{2}
\end{equation}
and that $k>0 $ is sufficiently large such that
\begin{equation}
\label{m1definition}
M_1 :=   \min \left\{\frac{3}{4} \, k\, a^2 - a - 1,\; k - 1 \right\} >0.
\end{equation}
Define
\begin{eqnarray}
\label{k1definition}
K_1  & := & \frac{1 + 2\,L^2}{M_1},
\\
\label{m0definition}
K_2 & :=  & \max \left\{ k\, a^2 + a + 1,\; k + 1 \right\} >0.
\end{eqnarray}
Then we have the inequalities
\begin{eqnarray}
\label{hilfe10082016oben}
M_1\,\int_0^L \left(u_t^2 +  u_x^2\right) \,dx
\leq  E_1(t)
& \leq &
K_2\,
\int_0^L \left(u_t^2 +  u_x^2\right) \,dx,
\\
\label{hilfe10082016}
\int_0^L \left(u_t^2+ (1+ 2\,L^2)\,  u_x^2\right) \,dx
 & \leq & \, K_1\,
E_1(t)
.
\end{eqnarray}
\end{lemma}
Note that inequality
(\ref{hilfe10082016})
is used in the proof of Theorem \ref{decay}
to obtain (\ref{27032017}).
{\bf Proof.}
Using Young's inequality we obtain
\begin{eqnarray*}
E_1(t) & \geq &
\int_0^L k \,\frac{3}{4} \, a^2 \, u_x^2 + k \, u_t^2
  - 2 \, \frac{a}{2}\, u_x^2 - u_x^2 - u_t^2  \, dx
  \\
  & = &\int_0^L \left( \frac{3}{4} \,k \, a^2   -   a  -  1 \right)\,u_x^2
    + \left( k- 1 \right) \,u_t^2  \, dx
    \\
    & \geq & M_1 \, \int_0^L u_x^2    + \,u_t^2  \, dx
    \\
    & \geq & \frac{M_1}{1 + 2L^2} \,\int_0^L \left(u_t^2+ (1+ 2\,L^2)\,  u_x^2\right) \,dx
    \end{eqnarray*}
    and the first inequality in (\ref{hilfe10082016oben}) and (\ref{hilfe10082016}) follow.
Using Young's inequality we also obtain
\begin{eqnarray*}
E_1(t) & \leq &
\int_0^L k \, a^2 \, u_x^2 + k \, u_t^2
  + 2 \, \frac{a}{2}\, u_x^2 + u_x^2  + u_t^2  \, dx
  \\
  & = &\int_0^L \left( k \, a^2    +  a  +1 \right)\,u_x^2
    + \left( k + 1\right) \,u_t^2  \, dx
    \\
    & \leq & K_2 \, \int_0^L u_x^2    + \,u_t^2  \, dx
    \end{eqnarray*}
    and the second inequality in (\ref{hilfe10082016oben}) follows.
  %

For $t\geq T_{period}>0$ we consider
the Lyapunov function
\begin{equation}
\label{edefinition}
E(t) = \int_{t - T_{period}}^t \,E_1(\tau)\, d\tau.
\end{equation}
By definition (\ref{edefinition}), $E(t)$ can be considered as a kind of
moving horizon time  average of $E_1(t)$.
We consider the time average since the
uncertain boundary data $b(\cdot, \omega)$
is distributed on the time
interval and condition (\ref{bassumption})
allows large values of
$|b(\cdot, \omega)|$ on  short time intervals.

We consider the situation where
on a finite time interval $[0,\, T]$,
the noise in the customer behavior approaches zero
with respect to the $H^1$-norm
as in (\ref{bassumption}).
We assume that during this process, the $C^2$--norm
$\|b(t,\,\omega)\|_{C^2(0,\, T)}$ is sufficiently small.
Then due to
 Theorem \ref{wellposedness},
a semi-global
$C^2$-solution exists
for sufficiently small initial data that are $C^2$--compatible
with the boundary conditions.
Moreover, the a priori estimate (\ref{apriori}) holds
which implies that by further decreasing the norms of
the initial and the boundary data, we can make the $C^2$--norm
of the solution $u$ as small as desired.

In Theorem \ref{decay} we state that
if $k>0$ is sufficiently large
under  appropriate  smallness conditions on $\bar u$ and $u$ and
if the uncertain customer profile
$b(\cdot,\,\omega)$ satisfies
(\ref{bassumption})
with sufficiently large $\nu$ and
$C_\nu$ sufficiently small,
our Lyapunov function
$E(t)$ defined in (\ref{edefinition})
decays exponentially fast
as in
(\ref{edecayint})
 with a
rate $\mu$ that is independent of $T$.

\begin{theorem}
\label{decay}
Let $L>0$ and $\lambda \in (\frac{1}{2},\,1)$ be given.
Assume that $k$ is sufficiently large such that
\begin{equation}
\label{sassumption}
k \geq \max
\left\{1,\;
\frac{4}{3} \left(\frac{1}{ a} +  \frac{1}{ a^2}\right)
,\;
\frac{1}{\lambda\, a}
\right\}
.
\end{equation}
Define
 \begin{equation}
\label{mudefinition}
\mu = \frac{1}{4\, {\rm e} \, L\, k}
\end{equation}
and the constant
\begin{equation}
\label{p0definition}
C_0
=
%
%
12\, k
+
4\,( k + 1)\,\left( 18 + 13 \,\theta +\frac{1}{a^2}\, ( 8  + 6\,\theta ) \right)
+ 10
.
\end{equation}
Let a  stationary state $\bar u(x)>0$,
$\bar u\in C^1(0,L)$ be given.
Assume that
\begin{equation}
\label{uquervoraus}
\bar u(x) \leq \min\left\{1,\;\frac{1}{4  \, k \,  {\rm e}}
,\; (1 -\lambda) \, \frac{a}{2}
,
\;
\frac{\mu}{C_0\, K_1}
\right\},\;
\bar u_x(x) \leq \min \left\{ 1,\; \frac{\mu}{C_0\, K_1}\right\}
\end{equation}
with  $K_1$ as defined in (\ref{k1definition}).
Let
$T_{period}>0$ and
$T>
T_{period}$ be given.
Assume that the initial data of system (\ref{sys}) and $b(\cdot,\,\omega)$
satisfy (\ref{epsilon0bedingung}) and the $C^2$--compatibility conditions such that
 Theorem \ref{wellposedness} implies that
 (\ref{sys})
 has a $C^2$--solution on $[0,\,T]\times [0,\, L]$ that satisfies the a-priori estimate
 (\ref{apriori}).
 Hence
we can assume that
 %
 $\|\varphi(x)\|_{C^2(0,\, L)}$,
 $\|\psi(x)\|_{C^1(0,\, L)}$
 and
 $\|b(t,\omega)\|_{C^2(0,\, T)}$
 are sufficiently small
  such that
\begin{equation}
\label{sassumptionb}
\left| u \right| \leq  \min\left\{ \bar u(0),\; (1 -\lambda) \, \frac{a}{2},\;
 \frac{1}{4  \, k \,  {\rm e}}
\right\},\;
\max\left\{ |u|, \,|u_x|,  \, |u_t|\right\}
\leq \min\left\{ 1,\;\frac{\mu}{C_0\, K_1}\right\}
.
\end{equation}
%
%

Assume that
 the uncertain function $b(t,\,\omega)$ satisfies
  (\ref{bassumption})
 with
 $\nu > \mu.$
 Define
 \begin{equation}
\label{deltadefinition}
 \delta = \nu - \mu >0.
\end{equation}
%
%
Then for all $t \in (
T_{period},\, T) $
for the solution of (\ref{sys}) we have the inequality
\begin{equation}
\label{edecayint}
E(t)
\leq
\exp\left(- \mu\; (t -  T_{period}) \right) \;
\left[ E(
 T_{period}) +
 \frac{C_g}{\delta}
 %
 \right]
%
\end{equation}
with $C_g>0$  defined by
\begin{equation}
\label{cgdefinition}
C_g
=
%
\left(
\frac{4}{3}\,{\rm e} \,a^2 \, k^2
+
\frac{1}{2 \,{\rm e}\, K_1\, k}
\right)
\,C_\nu
.
\end{equation}
In this sense $E(t)$ as defined in (\ref{edefinition})
decays exponentially with the rate
$\mu$
 that is independent of $T$.
%
%
%
If  we have
\begin{equation}
\label{nullvoraussetzung}
b(t,\,\omega)=0\; \mbox{\it for all }\;t \geq T- T_{period},
\end{equation}
then
\begin{equation}
\label{letzteungleichung}
\|u\|^2_{ H^1((T-T_{period},\, T)\times (0,\, L)) }
\leq
K_1 \;
\exp\left(- \mu\,
(T-T_{period}
)\right)
\left[
E(
T_{period})
+
 \frac{C_g}{\delta}
\right]
.
\end{equation}
%
%
\end{theorem}
\begin{remark}
Assumption (\ref{bassumption}) means
that the $H^1$--norm of the function that describes  the noise in the
customer behavior must  decay exponentially fast.
To guarantee the exponential decay of $E$ with the rate $\mu$
we assume that also the squared $H^1$--norm of $b(\cdot,\,\omega)$ decays
exponentially with a rate $\nu$  that is greater than ${\mu}$.
This condition
holds
if after
a finite time, 
the customer behavior becomes almost stationary.
\end{remark}
\begin{remark}
\label{remark2}
Since the conditions on $k$ do not depend on $T$,
the decay rate $\mu$ as defined in
(\ref{mudefinition}) does not depend on $T$.
Thus we can choose for example the time
$
T_{1/2} =
\frac{1}{\mu} \;{\rm ln}\left(
2\, K_1\, K_2
\right)
+
 T_{period}
.
$
Then
$\exp(-\mu \, T_{1/2}) = \frac{1}{2\, K_1\, K_2}
\exp(-\mu \, T_{period})
$.
For $t > T_{period} $,
we introduce the notation
 $X(t)=H^1((t -T_{period},\, t)\times (0,\, L))$.

If  for all $t \geq T- T_{period} $,
 condition (\ref{nullvoraussetzung}) holds
 and $T_{1/2} > T_{period}$,
 (\ref{letzteungleichung}) implies
 %
\[
\|u\|^2_{ X(T_{1/2}) }
%
\leq
\frac{1}{2}
\left[
\int_0^{T_{period}}
\int_0^L u_t^2(\tau,\, x) + u_x^2(\tau,\, x)\, dx \, d\tau
+
\frac{1}{K_2}\, \frac{C_g}{\delta}
\right].
\]
%
By a
trace
theorem
(see \cite{adams}, \cite{leoni})
we have the inequality
$
\int_0^L \left| u(t,\, x) \right|^2\, dx
\leq
C_e
\,
\|u\|_{ X(t)}
$
with an embedding constant
$C_e$.
Hence
(\ref{letzteungleichung})  yields an upper bound for
$\|u(T,\,\cdot)\|_{L^2(0,\, L)}$.
\end{remark}

\begin{remark}
The exponential decay
in  Theorem \ref{decay} can be interpreted as
an exponential decay of
the
function
\[
H(t)=
\int_{t-T_{period}}^t
\int_0^L u^2 + u_x^2 + u_t^2\, dx
\,
d\tau
=
\|u\|_{H^1((t-T_{period},\, t)\times (0,\, L))}^2.
\]
In fact, for all  $t\in (T_{period}, \,T)$
inequality (\ref{beweis17.08.16})
from the proof of Theorem \ref{decay} yields
\begin{equation}
\label{hilfsungleichung17.08.16}
H(t) \leq
K_1 \,
\exp\left(- \mu\; (t -  T_{period}) \right) \;
\left[ E( T_{period}) + \frac{C_g}{\delta}  \right]
 + 2 \, L\, C_\nu \, \exp(-\nu \, t)
 \end{equation}
 \[
=
 K_1 \,
\exp\left(- \mu\; (t -  T_{period}) \right) \;
\left[ E( T_{period}) + \frac{C_g}{\delta}
 + 2 \, \frac{L\, C_\nu}{K_1} \, \exp(- \delta \, t - \mu \,T_{period} ) \right].
 \]

On account of
the trace theorem mentioned in
Remark \ref{remark2}
and the definition of $H(t)$, (\ref{hilfsungleichung17.08.16}) implies
that  the $L^2$--norm of the
system state decays exponentially fast
with the rate $\mu$.
\end{remark}
\begin{remark}
\label{remark4}
It is interesting to compare the decay rate
$ \mu$ from (\ref{mudefinition}) with
the decay rate that is achieved for the system with the linear wave equation for $(t,\,x) \in (0,\infty)\times (0, L)$
that has already been studied in
\cite{ko:rapid}:
\[
\left\{
\begin{array}{l}
u(0,\, x) = \varphi(x)
\\
u_t(0,\, x) = \psi(x)
\\
u_{tt} - a^2\, u_{xx}=0
\\
u_x(t,\,0) = k \, u_t(t,\, 0)
\\
u(t,\, L)=0.
\end{array}
\right.
\]
The discussion in \cite{gubook}, (Chapter 5.2)
implies
that
for $k > \frac{1}{a}$,
$\varphi \in H^1((0,\,L))$
with $\varphi(L)=0$
and $\psi \in L^2((0,\, L))$
the classical energy $E_{classic}$
decays exponentially and the optimal rate is
$
\mu_0 =
\frac{a}{L} \,
{\rm ln}\left(
1  + \frac{2}{a\, k - 1}
 \right)
.
$
Thus we have
\[\frac{\mu_0}
{\mu}
=
4\, {\rm e} \,
{\rm ln}\left(
\left(1 + \frac{2}{a\, k - 1}\right)^{ak} \right).
\]
Hence
$\lim_{(a\, k) \rightarrow \infty}
\frac{\mu_0}
{\mu}
= 8 \, e.
$
%
Thus asymptotically for large values of $(a\, k)$
the decay rate $\mu$ for the quasilinear system
differs from the rate for the linear wave equation
only by a multiplicative constant
that is close to ${8\,{\rm e}}$.

\end{remark}

For the proof of Theorem \ref{decay},
we use the following variant of
Gronwall's Lemma, that we prove for the convenience of the reader.

\begin{lemma}[Gronwall's Lemma]
\label{lemma6}
Let  real numbers $\mu>0$,
$\nu>0$ and $C_g>0$
be given such that
$ \nu > \mu$.
Define
\[\delta = \nu - \mu >0.\]
Assume that $U(t)\geq 0$ is a differentiable
function that satisfies for all
 $t\in [0,\, T]$ the inequality
\begin{equation}
\label{diffinequality}
U'(t) \leq - \mu \, U(t) + C_g  \exp(- \nu \, t).
\end{equation}

Then for all  $t\in [0,\, T]$  the function $U(t)$
satisfies the inequality
\begin{equation}
\label{gronwallungleichung}
0 \leq U(t) \leq \exp(- \mu\, t) \left( U(0) + \frac{C_g}{\delta}\right).
\end{equation}
\end{lemma}
%
{\bf Proof. }
We define the auxiliary function
\begin{equation}
{\cal H}(t) = \exp( \mu\, t)\; U(t).
\end{equation}
Then we have ${\cal H}(0)= U(0)$.
The product rule and (\ref{diffinequality}) imply
\begin{eqnarray*}
{\cal  H}'(t) & = & \mu\, {\cal  H}(t) + \exp( \mu\, t)\; U'(t)
\\
  & \leq  & \mu \, {\cal H}(t) + \exp( \mu\, t)\; \left[   - \mu \, U(t) + C_g  \exp(-  \nu \, t) \right]
\\
  & = & \mu \, {\cal H}(t)    - \mu \, {\cal H}(t)  + C_g  \exp( (\mu -  \nu) \, t)
\\
 & = &  C_g \, \exp( - \delta \, t)  .
\end{eqnarray*}
By integration the inequality
${\cal H}'(t) \leq  C_g  \exp( - \delta \, t) $ yields
\begin{eqnarray*}
{\cal H}(t) - {\cal H}(0) & = &  \int_0^t {\cal H}'(\tau) \, d\tau
\\
&\leq &
\int_0^t   C_g  \,\exp( - \delta \, \tau)\,d\tau
=
   C_g  \,   \frac{1}{\delta} \, ( 1 - {\rm e}^{- \delta\, t}).
\end{eqnarray*}
Hence we have
\begin{eqnarray*}
U(t)  & = & {\rm e}^{-\mu \, t}  \,{\cal H}(t)
\\
&\leq &
{\rm e}^{- \mu \,t} \left( {\cal H}(0) +   C_g  \,   \frac{1}{\delta} ( 1 - {\rm e}^{- \delta\, t})  \right)
\\
 & \leq &
{\rm e}^{- \mu \,t} \left( U(0) +   \frac{ C_g  }{\delta}   \right)
\end{eqnarray*}
and (\ref{gronwallungleichung}) follows.
%

In the subsequent analysis an upper bound for $F$ is used that is presented in the following lemma.
\begin{lemma}
\label{flemma}
Define
\begin{equation}
\label{tlidefinition}
T_{Li}(t) = \max_{x\in [0,L]} \left\{ |u(t,x)|, |u_x(t,x)|,  |u_t(t,x)|,\, |\bar u(x)|,\, |\bar u_x(x)|\right\}.
\end{equation}
Assume that
\begin{equation}
\label{tliassumption}
\bar u + u \geq 0,\; 0 < \bar u \leq \frac{a}{2} ,\; T_{Li}(t)\leq 1.
\end{equation}
Then we have the equation
\begin{equation}
\label{fgleichung}
F   =  - 2 \, u_t\, (u_x + \bar u_x) - \theta\, (u + \bar u)\, u_t
- 2 \, u \, (u_x + \bar u_x)^2 - 4\, \bar u \, \bar u_x \, u_x - 2 \, \bar u \, u_x^2
\end{equation}
\begin{eqnarray*}
& -
&  \frac{3}{2}\, \theta \, u\, ( u + 2\, \bar u) \, (u_x + \bar u_x)
 -
\frac{3}{2} \theta  \, \bar u^2
\, u_x
-
\frac{2\, u \, \bar u + u^2}{a^2 - \bar u^2} \left(
2 \, \bar u \, \bar u_x^2 + \frac{3}{2}\, \theta \, \bar u^2
\, \bar u_x \right).
\end{eqnarray*}
and the upper bound
\begin{equation}
\label{fbound}
\left|F(x,u(t,\cdot),u_x(t,\cdot),u_t(t,\cdot))\right|
\end{equation}
\[\leq
\left[ 18 + 13 \,\theta +\frac{1}{a^2}\, ( 8  + 6\,\theta ) \right]
\;
 T_{Li}(t)
\;
 \left( |u(t,x)| + |u_x(t,x)| +  |u_t(t,x)|\right).
\]
\end{lemma}
{\bf Proof. }
 Equation (\ref{fgleichung}) follows with the definition of $\tilde F$,
using the assumptions $\bar u + u >0$, $\bar u>0$
 from (\ref{tliassumption})
to eliminate the absolute value brackets.
From (\ref{fgleichung})  we obtain the upper bound (\ref{fbound})
using again the conditions from (\ref{tliassumption}) that imply
$T_{Li}(t)^2\leq T_{Li}(t)$.
%

Now we can proceed with the proof of
 Theorem \ref{decay}, where we show that $E(t)$ satisfies
 a differential inequality of the type (\ref{diffinequality})
 and thus
 due to
 (\ref{deltadefinition})
 we can apply Lemma \ref{lemma6}.

{\bf Proof of Theorem \ref{decay}.}
\label{lemmaE1}
Define the exponential weight
\[h_2(x)= \exp\left(-\frac{x}{L}\right),\; x\in [0,\, L].\]
The time-derivative of $E_1$ is given by the
 equation
\begin{equation}
\label{e1strichgleichung}
\frac{d}{dt}E_1(t)
 =  I_1 + I_2 + I_3
\end{equation}
with
\begin{eqnarray}
\label{5.15}
I_1 &= &
\int_0^L h_{2x} \, (a^2-(\bar u+u)^2) \,  u_x^2  +
 h_{2x} \, u_t^2 \,dx,
\\
\label{5.16}
I_2& = &
\int_0^L
2 \, k\, (\bar u_x+u_x) \, u_t^2
-2\, k\, (\bar u+u)\, u_t \, u_x^2+ 4 \, k\, (\bar u+u)(\bar u_x+u_x) \, u_x  \, u_t
\\
\nonumber
& + & 2 \, k  \, F \, u_t
 -2 \, h_2 \, u_t \, u_x^2 - 2 \, h_2 \, (\bar u+u)(\bar u_x+u_x) u_x^2-2 \, h_2 \, F \, u_x \,dx,
 \\
\label{5.17}
I_3 &= &
[(a^2-(\bar u+u)^2)(2 \,
\,k \, u_x \, u_t-h_2 \, u_x^2)-(2\,
\,k\, (\bar u+u)+h_2) \,u_t^2]_{x=0}^L.
\end{eqnarray}
This can be seen as follows.
%
%
With the notation
\begin{equation}
\label{dhutdefinition}
\hat d = a^2-(\bar u+u)^2
\end{equation}
we have
$\hat d_t = - 2 \, (\bar u+u) \, u_t$,
$\hat d_x = - 2 \, (\bar u+  u ) \,(\bar u_x+  u_x ) $
 and
\[
E_1(t)=
\int_0^L
k\, \Big(  \hat d\, u_x^2 +  u_t^2\Big)-2 \,  h_2 \,
\Big((\bar u+u)\, u_x^2+u_t \, u_x\Big)\,dx.
\]
Hence differentiation yields
\begin{eqnarray*}
\frac{d}{dt}E_1(t) & = &
\int_0^L
2\,k \,  \Big[
(u_{tt}
-
(\bar u + u) \, u_x^2 )
\, u_t +
\hat d\,  u_x\, u_{xt}
\Big]
\\
&
-
& 2\,  h_2 \,
\Big[u_t \, u_x^2 + (  u_{tt} + 2 (\bar u+u)\, u_{xt}) \, u_x
 +   u_t\, u_{xt}  \Big]\,dx.
\end{eqnarray*}
Now integration  by parts
for the term  $\hat d\,  u_x\, u_{xt} = \left(\hat d\,  u_x \right)\, \left(u_t\right)_x$
yields the equation
\begin{eqnarray*}
\frac{d}{dt}E_1(t)
& = &
\int_0^L
2\,k \, \Big[u_{tt}
-
\hat d\,  u_{xx}
-
\hat d_x\,  u_x\,
-
(\bar u + u) \, u_x^2
\Big]\, u_t
\\
&
-
&
2\,  h_2 \,
\Big[u_t \, u_x^2 + (  u_{tt} + 2 (\bar u+u)\, u_{xt}) \, u_x
 +   u_t\, u_{xt}  \Big]\,dx
 +
\left[ 2\,k\, \hat d\,  u_x\, u_{t}\right]_{x=0}^L.
\end{eqnarray*}
Hence we get  the equation
\begin{eqnarray*}
\frac{d}{dt}E_1(t)
& =  &\int_0^L 2 \, k
\left[(u_{tt}- \hat d\, u_{xx})\, u_t
-\hat d_x
\, u_x  \, u_t
-(\bar u+u)  u_x^2 u_t \right]
\\
&
-
&
2 \, h_2 \,
\left[ u_t \, u_x^2+(u_{tt}+2\,
(\bar u+u)\, u_{tx}) \, u_x + u_t \, u_{tx}
\right]
\, dx
\\
&
+
&
\left[2 \, k \, \hat d \, u_x \, u_t\right]_{x=0}^L.
\end{eqnarray*}
By the partial differential equation in (\ref{sys})
we
have
$u_{tt}- \hat d\, u_{xx}
=
F-2 \, (\bar u+u) \,u_{tx}
$
and
obtain
\begin{eqnarray*}
\label{5.14}
\frac{d}{dt} E_1(t)
& =  & \int_0^L 2 \, k\,
\left[(F-2\, (\bar u+u) \, u_{tx}) \, u_t+2\,
(\bar u+u)\,(\bar u_x+u_x) \, u_t \, u_x - (\bar u+u) u_t \, u_x^2
\right]\nonumber\\
& -
&2 \, h_2
\left[
(F+
\hat d\,
u_{xx})\,u_x + u_t \, u_x^2 + u_t  \, u_{tx}
\right] \, dx
+
\left[
2 \, k \,
\hat d
\,
 u_x
 \,
  u_t
\right]_{x=0}^L.
\nonumber
\end{eqnarray*}
Using integration by parts we obtain the identities
\begin{eqnarray*}
\int_0^L -4\, k \, (\bar u+u)u_t u_{tx} \, dx & = &
\int_0^L - 2 \, k\,(\bar u+u) (u_t^2)_x \, dx
\\
& = & [-2\, k\, u_t^2(\bar u+u)]_{x=0}^L +\int_0^L 2\, k \, (\bar
u+u)_x  \, u_t^2\, dx
\end{eqnarray*}
and
\begin{eqnarray*}
\int_0^L -2 \, h_2 \, \hat d\,u_x u_{xx} - 2 \, h_2 \, u_t \, u_{tx} \, dx  & = &
\int_0^L - h_2   \, \hat d\, (u_x^2)_x  - h_2 \, (u_t)^2_x\, dx
\end{eqnarray*}
\begin{eqnarray*}
& = &
[-h_2  \, \hat d\, u_x^2  - h_2 \, u_t^2]_{x=0}^L
\\
&  + &
\int_0^L h_{2x}   \, \hat d\, u_x^2
- 2 \, h_2(\bar
u+u)(\bar u_x+u_x) \, u_x^2 + h_{2x} \, (u_t)^2\, dx.
\end{eqnarray*}
Using these identities we obtain equation (\ref{e1strichgleichung}).
Here, $I_3$ contains all the terms coming from the boundary
and $I_1= \int_0^L h_{2x} \,\hat d\, u_x^2 + h_{2x} \,u_t^2 \, dx$
contains all the terms where $h_{2x}$ appears.
The remaining terms appear in $I_2$.


We use the notation $\bar u_0  = \bar u(0)$ and $\bar u_L  = \bar u(L)$.
We have
$I_3 = I_3^L - I^0_3$
with
\begin{equation}
\label{i30defi}
I_3^0  =   \Big( a^2 - \left(  (\bar u_0 + u(t,0) )+ \frac{1}{k} \right)^2 \Big) \, u_x^2(t,0)
\end{equation}
and
\begin{eqnarray*}
I_3^L
& = & \left[a^2-(\bar u_L+u(t,\, L))^2\right]\,\left[2 \, k \, u_x(t,\, L) \, u_t(t,\, L) - \exp(-1) \, u_x(t,L)^2\right]
\\
& - & \left[2\, k\, (\bar u_L + u(t,\, L) )+ \exp(-1) \right] \;u_t(t,\, L)^2.
\end{eqnarray*}
Due to
(\ref{sassumption})
 we have
 \[\frac{1}{k} \leq \lambda\, a.\]
 Moreover, due to
 (\ref{uquervoraus}) and (\ref{sassumptionb})
 we have
 \begin{equation}
 \label{ahalbeungleichung}
 |\bar u + u|\leq \bar u + |u|\leq
 (1-\lambda)\, a
 \leq \frac{a}{2}.
 \end{equation}
Hence we have
\[
\left( \frac{1}{k} + (\bar u + u)\right)^2 \leq a^2.
\]
 %
  Thus we have
 \begin{equation}
 \label{i30ungleichung}
 I_3^0\geq 0.
 \end{equation}
%
%
%
%
%
Moreover, due to
(\ref{ahalbeungleichung}),  the Cauchy--Schwarz inequality
and
 (\ref{bassumption})
we have
\begin{eqnarray*}
& &
\int_{t-T_{period}}^t (a^2-(\bar u_L+u(\tau,\, L))^2)\; 2 \, k \, u_x(\tau,\, L) \, b_t(\tau,\, \omega)\, d\tau
\\
&
\leq
&
2\; a^2  \, k \,\int_{t - T_{period}}^t  \left|u_x(\tau,\, L)\right|\, \,  \left|b_t(\tau,\, \omega) \right|\, d\tau
\\
&
\leq
&
2\; a^2  \, k \,\sqrt{\int_{t -T_{period}}^t  \left(u_x(\tau,\, L)\right)^2
\,d\tau}
\,
\sqrt{
\int_{t-T_{period}}^t
 \, \left(b_t(\tau,\, \omega)\right)^2 \, d\tau}
 \\
 & \leq &
 2\; a^2  \, k\;
\sqrt{C_\nu}
\;
\sqrt{\int_{t -T_{period}}^t  \left(u_x(\tau,\, L)\right)^2
\,d\tau}
\exp\left( - \frac{\nu}{2} t
\right)
.
\end{eqnarray*}
Due to (\ref{ahalbeungleichung}) we have
\begin{equation}
\label{dreiviertelungleichung}
 \frac{3}{4} \, a^2 \leq a^2-(\bar u+u)^2\leq  a^2.
 \end{equation}
Hence
we obtain
\begin{equation}
\label{ungleichung8916}
\int_{t-T_{period}}^t
I_3^L \, d\tau
\end{equation}
\[
\leq
 2\; a^2  \, k\;\sqrt{C_\nu}
 \frac{
\sqrt{\int_{t -T_{period}}^t  \left(u_x(\tau,\, L)\right)^2
\,d\tau}
}{
\exp\left( \frac{\nu}{2} \, t
\right)}
-
\frac{3\,a^2}{ 4\, {\rm e}}
\int_{t -T_{period}}^t  \left(u_x(\tau,\, L)\right)^2
\,d\tau
\]
\[
=
\frac{2\; a^2  \, k\;\sqrt{C_\nu}}
{\exp\left(\frac{\nu}{2} \, t \right)}
\,
 \|u_x(\cdot,\, L)\|_{L^2(t -T_{period},\, t)}
-
 \frac{3\,a^2}{ 4\, {\rm e}}\,
\|u_x(\cdot,\, L)\|_{L^2(t -T_{period},\, t)}^2
.
\]
Define the polynomial
\begin{equation}
\label{pdefinition}
p_3(z) = 2\; a^2  \, k\;\sqrt{C_\nu}\, \exp\left( -\frac{\nu}{2} \, t \right)
z -
\frac{3\,a^2}{ 4\, {\rm e}}\,z^2.
\end{equation}
The polynomial $p_3$ has the form
$p_3(z)= \alpha \, z - \beta \, z^2$ with $\beta>0$.
Hence for all $z\in (-\infty,\,\infty)$, we have the inequality
\begin{equation}
\label{pungleichung}
p_3(z)\leq  p_3(\frac{\alpha}{2 \, \beta}) = \frac{\alpha^2}{4\, \beta}
=\frac{4}{3}\,{\rm e} \,a^2 \, k^2 \,C_\nu\, \exp(-\nu \, t).
\end{equation}
Hence (\ref{ungleichung8916}) implies
\begin{equation}
\label{immerungleichung}
\int_{t - T_{period}}^t I_3^L \, d\tau \leq
\frac{4}{3}\,{\rm e} \,a^2 \, k^2 \,C_\nu\, \exp(-\nu \, t)
.
\end{equation}
%
Thus  due to (\ref{i30ungleichung})  we have
\begin{equation}
\label{i3ungleichung}
\int_{t-T_{period}}^t I_3 \, d\tau
 \leq
 \frac{4}{3}\,{\rm e} \,a^2 \, k^2 \,C_\nu\, \exp(-\nu \, t)
 .
\end{equation}
We have
\begin{equation}
\label{i1ungleichung9816}
I_1\leq -  \frac{1}{2\,{\rm e} \, L\,  k }\, E_1(t).
\end{equation}
This can be seen as follows. We have
\begin{eqnarray*}
I_1  & = &
-\frac{1}{L}
\int_0^L \exp\left(-\frac{x}{L}\right) \, (a^2-(\bar u+u)^2) \,  u_x^2  +
\exp\left(-\frac{x}{L}\right) \, u_t^2 \,dx
\\
& \leq &
-\frac{1}{2 \, k \,L\,{\rm e}}
\int_0^L  k [ (a^2-(\bar u+u)^2) \,  u_x^2  +   u_t^2 ]\,dx
\\
& - &
\frac{1}{2 \,L}
\int_0^L \exp\left(-\frac{x}{L}\right) \, (a^2-(\bar u+u)^2) \,  u_x^2  +
\exp(-\frac{x}{L}) \, u_t^2 \,dx.
\end{eqnarray*}
Since $\bar u>0$ is strictly increasing and due to (\ref{sassumptionb}) we have
\begin{equation}
\label{groessergleichnullungleichung}
\bar u(x) + u(x) \geq \bar u(0) - |u(x)| \geq
\bar u(0) - \bar u(0) =0,
\end{equation}
 hence
\begin{eqnarray*}
I_1
& \leq &
-\frac{1}{2\,{\rm e}\,  k \,L}
\int_0^L  k [ (a^2-(\bar u+u)^2) \,  u_x^2  +   u_t^2 ]
-2
\exp\left(- \frac{x}{L}\right)\,
\left[
(\bar u+u)\,u_x^2
+ u_t \, u_x\right]
\,dx
\\
& + &
\frac{1}{2\, {\rm e}\,  k \,L}
\int_0^L
\exp\left(- \frac{x}{L}\right)
\left[
(u_t^2 +  u_x^2)
- k \,{\rm e}\, [ (a^2-(\bar u+u)^2) \,  u_x^2  +   u_t^2 ]
\right]
\, dx.
\end{eqnarray*}
Due to  (\ref{sassumption}), we have $ 1 - k\,{\rm e} < 0$, hence
\[
I_1
 \leq
-\frac{1}{2\,{\rm e}\,  k \,L} \, E_1(t)
 +
\frac{1}{2\, {\rm e}\,  k \,L}
\int_0^L
\exp\left(- \frac{x}{L}\right)
\left[
 1
- k \,{\rm e}\,  (a^2-(\bar u+u)^2)
\right]\,  u_x^2
\, dx.
\]
Due to
(\ref{dreiviertelungleichung})
and  (\ref{sassumption}) this yields
\begin{eqnarray*}
I_1
& \leq &
-\frac{1}{2\,{\rm e}\,  k \,L} \, E_1(t)
 +
\frac{1}{2\, {\rm e}\,  k \,L}
\int_0^L
\exp\left(- \frac{x}{L}\right)
\left[
 1
- k \,{\rm e}\,  \frac{3}{4} \, a^2
\right]\,  u_x^2
\, dx
\\
& \leq &
-\frac{1}{2\,{\rm e}\,  k \,L} \, E_1(t).
\end{eqnarray*}
Thus we have shown (\ref{i1ungleichung9816}).

From (\ref{uquervoraus}) we have $\bar u \in (0, \, a/2)$ and
(\ref{sassumptionb}) implies that
\[\max\{|\bar u|,\, \, |\bar u_x|\} \leq \min\left\{ 1,\; \frac{\mu}{K_1\, C_0}\right\}\]
with $C_0$ as defined in (\ref{p0definition}).
Moreover  (\ref{sassumptionb}) implies
\[\max\{|u|,\, \, | u_x|,\, |u_t|\} \leq \min\left\{ 1,\; \frac{\mu}{K_1\, C_0}\right\}.\]
Hence for $T_{Li}(t)$ as defined in
(\ref{tlidefinition}) we have
\begin{equation}
\label{tliungleichung}
T_{Li}(t) \leq  \min\left\{ 1,\; \frac{\mu}{K_1\, C_0}\right\}.
\end{equation}
Hence (\ref{groessergleichnullungleichung}) implies that (\ref{tliassumption}) holds.
Thus we can use (\ref{fbound}) to derive an upper bound for $I_2$.
Due to (\ref{tliungleichung})  we have
\begin{eqnarray*}
I_2& = &
\int_0^L
2 \, k\, (\bar u_x+u_x) \, u_t^2
-2\, k\, (\bar u+u)\, u_t \, u_x^2+ 4 \, k\, (\bar u+u)(\bar u_x+u_x) \, u_x  \, u_t
\nonumber
\\
& + & 2 \, k  \, F \, u_t
 -2 \, h_2 \, u_t \, u_x^2 - 2 \, h_2 \, (\bar u+u)(\bar u_x+u_x) u_x^2-2 \, h_2 \, F \, u_x \,dx
 \\
 & \leq  &
\int_0^L
4 \, k \, T_{Li}(t)\, u_t^2
 + 4 \, k\,  T_{Li}(t)\, u_x^2
 + 8\, k\,  T_{Li}(t)\,(u_x^2 +  u_t^2) + 2 \,T_{Li}(t) \, u_x^2 +  8 \, T_{Li}(t)\, u_x^2
\nonumber
\\
& + & 2 \, k  \, F \, u_t
 -2 \, h_2 \, F \, u_x \,dx
  \\
 & \leq  &
\int_0^L \left( 12 \, k + 10 \right) \,  T_{Li}(t)\,  \left(  u_t^2 + u_x^2\right)
+
 2 \, k  \, F \, u_t
 -2 \, h_2 \, F \, u_x \,dx.
\end{eqnarray*}
Hence (\ref{fbound}) implies
%
\begin{eqnarray*}
I_2 & \leq  &
\int_0^L \left( 12 \, k + 10 \right) \,   T_{Li}(t)\, \left(  u_t^2 + u_x^2\right)
+
\end{eqnarray*}
\[+
2\, k \,
\left[ 18 + 13 \,\theta +\frac{1}{a^2}\, ( 8  + 6\,\theta ) \right]
\;
 T_{Li}(t)
\;
 \left( |u\, u_t| + |u_x\, u_t| +  |u_t^2|\right)
\]
\[
+
2\,
\left[ 18 + 13 \,\theta +\frac{1}{a^2}\, ( 8  + 6\,\theta ) \right]
\;
 T_{Li}(t)
\;
 \left( |u\, u_x| + |u_x^2| +  |u_t\,u_x|\right)\,dx.
\]
Using Young's inequality, this yields with the definition of $C_0$ in (\ref{p0definition})
\begin{eqnarray*}
I_2 &  &
\end{eqnarray*}
\[
\leq
\int_0^L \left[
12 \, k + 10
+
4 \left(  k + 1\right)\,
\left(
18 + 13 \,\theta + \frac{1}{a^2}\,( 8 + 6 \, \theta)
\right)
\right] \,
 T_{Li}(t)\,
 \left(u^2 +   u_t^2 + u_x^2\right)\, dx
\]
\begin{equation}
\label{18.08.2016}
\leq
C_0\;  T_{Li}(t)\,\int_0^L
 \left(u^2 +   u_t^2 + u_x^2\right)\, dx.
\end{equation}
%
We have
\[u(t,\, x) = u(t,\, L) -\int_x^L u_x(t,\, x)\, dx = b(t,\,\omega)  -\int_x^L u_x(t,\, x)\, dx\]
hence we have
\[\left|u(t,\, x)\right|
 \leq \left| b(t,\,\omega)\right|  + \sqrt{L}\,  \left(\int_0^L |u_x(t,\, x)|^2\, dx\right)^{1/2}.
\]
This implies with Young's inequality
\[
\left|u(t,\, x)\right|^2 \leq
 2 \left| b(t,\,\omega)\right|^2  +2\,L\, \int_0^L |u_x(t,\, x)|^2\, dx.
\]
Thus we obtain
\begin{equation}
\label{poincare}
\int_0^L \, \left|u(t,\, x)\right|^2 \, dx \leq
  2\, L\,  \left| b(t,\,\omega)\right|^2  +2\,L^2\, \int_0^L |u_x(t,\, x)|^2\, dx.
\end{equation}
Hence (\ref{18.08.2016}) implies  
\[
I_2
\leq
 C_0 \,  T_{Li}(t) \,
 \left[ \int_0^L \left(u_t^2+ (1+ 2\,L^2)\,  u_x^2\right) \,dx
 + 2\, L\,  \left| b(t,\,\omega)\right|^2\right].
\]
Due to (\ref{groessergleichnullungleichung}) and (\ref{ahalbeungleichung}),
(\ref{110816}) holds.
Condition (\ref{sassumption}) implies that (\ref{m1definition}) is valid.
Thus we can apply (\ref{hilfe10082016})  to obtain
\begin{equation}
\label{27032017}
I_2
\leq
C_0\, T_{Li}(t) \,\left[ \;K_1 \, E_1(t) + 2\, L\,  \left| b(t,\,\omega)\right|^2 \right]
\end{equation}
with the constant $K_1$ from Lemma \ref{fragezeichen}.
Thus  (\ref{i1ungleichung9816}) and
(\ref{tliungleichung})
imply
\[I_1 + I_2 \leq
-\frac{1}{2\,{\rm e}\,  k \,L} \, E_1(t) +
\mu \, E_1(t)
+
\frac{1}{2\,{\rm e} \,K_1\,  k} \, \left| b(t,\,\omega)\right|^2
\]
\[
=
-  \mu\, E_1(t)
+
\frac{1}{2\,{\rm e}\,K_1\,  k} \, \left| b(t,\,\omega)\right|^2.
\]
Thus we have
\[\frac{d}{dt}E_1(\tau)
\leq
-  \mu\, E_1(\tau) + I_3 +
\frac{1}{2\,{\rm e}\,K_1\,  k} \, \left| b(t,\,\omega)\right|^2.
\]
This implies in turn  due to
(\ref{i3ungleichung}) and (\ref{bassumption})
\begin{eqnarray*}
E_1(t ) -  E_1(t-  T_{period})
& = &
\int_{t - T_{period}}^t
\frac{d}{dt}E_1(\tau)
\, d\tau
\end{eqnarray*}
\[
 \leq
-  \mu\,
\int_{t - T_{period}}^t
 E_1(\tau)
\,d\tau
 +
\left[
\frac{4}{3}\,{\rm e} \,a^2 \, k^2 \,C_\nu
+
\frac{C_\nu}{2 \,{\rm e}\, K_1\, k} \,
\right]
\,
\exp\left( - {\nu} \,t
\right).
\]
Thus for $E(t)$ as defined in (\ref{edefinition})
for all $t\geq T_{period}$ we have the inequality
\[\frac{d}{dt}
E(t) =
E_1(t) - E_1(t -T_{period})
\leq
-  \mu\,
E(t)
+
C_g\; \exp\left( - {\nu} \, t\right)
\]
with the constant $C_g>0$ as defined in (\ref{cgdefinition}).

For $s\geq 0$, let $U(s) = E (T_{period} + s)$.
By  Lemma \ref{lemma6}, this implies
(\ref{edecayint}).

Due to
(\ref{poincare}) and
(\ref{hilfe10082016}),
inequality  (\ref{edecayint}) yields for $t\geq T_{period}$
\[
\int_{t-T_{period}}^t
\|u(\tau,\cdot)\|_{H^1(0, \, L)}^2
+
\|u_t(\tau,\cdot)\|_{L^2(0, \, L)}^2
\,d\tau
\]
\[
\leq
\int_{t-T_{period}}^t
\int_0^L ( 1 + 2\, L^2) \, u_x^2(\tau,\, x) + u_t^2(\tau,\, x) dx
+
2\, L\, | b(\tau,\,\omega)|^2
 \, d\tau
\]
\[
\leq
K_1 \;
\int_{t-T_{period}}^t \,
E_1(\tau)\, d\tau
+
2\, L\, C_\nu\, \exp(-\nu \,t)
\]
\begin{equation}
\label{beweis17.08.16}
\leq
K_1 \;
\exp\left(- \mu\,(t-T_{period})\right)
\left(E(T_{period}) + \frac{C_g}{\delta} \right)
+2\, L\, C_\nu\, \exp(-\nu \,t).
\end{equation}
Analogously, (\ref{nullvoraussetzung}) yields
%
\[
\int_{T-T_{period}}^T
\|u(\tau,\cdot)\|_{H^1(0, \, L)}^2
+
\|u_t(\tau,\cdot)\|_{L^2(0, \, L)}^2
\,d\tau
\]
\[
\leq
K_1 \;
\exp\left(- \mu\,
(T-T_{period}
)\right)
\left(E(T_{period}) + \frac{C_g}{\delta} \right)
\]
and hence (\ref{letzteungleichung}).
%
Thus we have proved Theorem \ref{decay}. 

%

\section{Conclusions}
In many applications,
uncertainty in the boundary data occurs, for example if
uncertain customer behavior influences
the boundary data of  the system.
We have shown that also with such  uncertain boundary data
on parts of the boundary,
a boundary--feedback law can lead to exponential
decay  for a system that
is governed by a quasilinear hyperbolic equation,
provided that the perturbations of the boundary data
decay exponentially  with a decay rate that is sufficiently large.
To guarantee the exponential decay,
the feedback parameter has to be sufficiently large
and the desired stationary state has to be sufficiently small.
Due to the nonlinearity of the system, the result is local,
that is the initial state has to be sufficiently close to the desired
stationary state.
The proof is based upon a suitably defined Lyapunov function.
Our results  show that with the feedback controller the energy
of the system in
the sense of our Lyapunov function decays exponentially
fast, even if there is some unknown input at parts of the boundary.
The exponential decay of the
Lyapunov function also implies that the $L^2$--norm of the
system state decays exponentially fast.

If the boundary perturbations do not decrease exponentially fast to zero,
is is not possible to achieve exponential decay of the system state.
However, we expect that with the feedback law
(\ref{0rb}),
also in this case
the system state will decay
with the same speed as the
boundary perturbations.
A detailed study of this situation is
a project for future research. In such a situation
it could be useful to
consider  the feedback law
(\ref{0rb})  in the integral form
\[\tilde u(t,0) =
\tilde u(0,\,0)  +
\frac{1}{k} \int_0^t \tilde u_x(\tau,\, 0) - \bar u_x( 0) \, d\tau.
\]
The Riemann invariants of (\ref{2.1}), (\ref{2.2})
are
$R_\pm = - \frac{q}{\rho} \mp a\,{\rm ln}(\rho)$,
see \cite{dick}.
For the velocity this yields $\tilde u = - \frac{1}{2}(R_+ + R_-)$.
Hence the feedback law (\ref{0rb})
is equivalent to the linear Riemann feedback
\[
(R_+)_x(t,\,0) - k (R_+)_t(t,\,0) = - 2 \, \bar u_x(0)
 - (R_-)_x(t,\,0) + k (R_-)_t(t,\,0).
\]
Therefore we hope that our analysis
is a motivation to consider this type of Riemann
feedback laws in future studies.

\section*{Acknowledgments}
We want to thank the referees for their valuable comments.
This work is  supported by DFG in the framework of the Collaborative Research Centre
CRC/Transregio 154, Mathematical Modelling, Simulation and Optimization Using the Example of Gas Networks, project C03.


\end{document}